\documentclass[notitlepage]{amsart}

\usepackage{amsfonts}
\usepackage[latin1]{inputenc}
\usepackage[all]{xy} 
\usepackage{latexsym,amssymb,amscd}
\usepackage{amsmath}
\newtheorem{pro}{Proposition}[section]
\newtheorem{teo}[pro]{Theorem}
\newtheorem{defi}[pro]{Definition}
\newtheorem{lem}[pro]{Lemma}
\newtheorem{cor}[pro]{Corollary}
\newtheorem{rk}[pro]{Remark}

\newcommand{\Hom}{\mathrm{Hom}}
\newcommand{\thick}{\overline{\Delta}_{\mathcal{T}}}
\newcommand{\thickH}{\overline{\Delta}_{\mathbf{D}^b(\mathcal{H})}}

\newcommand{\la}{\left\langle}
\newcommand{\ra}{\right\rangle}
\newcommand{\U}{\mathcal{U}}

\newcommand{\sU}{\overline{\mathcal{U}}}
\newcommand{\T}{\mathcal{T}}

\newcommand{\HH}{\mathcal{H}}
\newcommand{\C}{\bf{C}}

\newcommand{\A}{\mathcal{A}}
\newcommand{\B}{\mathcal{B}}

\newcommand{\X}{\mathcal{X}}
\newcommand{\Y}{\mathcal{Y}}

\newcommand{\N}{\mathbb{N}}
\newcommand{\Enteros}{\mathbb{Z}}
\newcommand{\pd}{\mathrm{pd}}
\newcommand{\id}{\mathrm{id}}
\newcommand{\resdim}{\mathrm{resdim}}
\newcommand{\coresdim}{\mathrm{coresdim}}
\newcommand{\add}{\mathrm{add}}
\newcommand{\mini}{\mathrm{min}}
\newcommand{\maxi}{\mathrm{max}}
\newcommand{\D}{\mathbf{D}^b}
\newcommand{\Kb}{\mathbf{D}^b\,(\mathrm{mod}\,(\Lambda))}
\newenvironment{dem}{\noindent\bf Proof. \rm }{$\ \Box$}
\usepackage{latexsym,amssymb,amscd}
\usepackage{amsmath}

\begin{document}
\title[Auslander-Buchweitz context and co-$t$-structures]{Auslander-Buchweitz context and co-$t$-structures}
\author{O. Mendoza, E. C. S\'aenz, V. Santiago, M. J. Souto Salorio.}
\thanks{
The authors thank the financial support received from
Project PAPIIT-UNAM IN101607 and MICINN-FEDER   TIN2010-18552-C03-02.}
\date{}
\begin{abstract}
We show that the relative Auslander-Buchweitz context on a triangulated category $\T$ coincides with the notion of
co-$t$-structure on certain triangulated subcategory of $\T$ (see Theorem \ref{M2}). In the Krull-Schmidt case, we stablish a bijective correspondence between co-$t$-structures and cosuspended, precovering subcategories (see Theorem \ref{correspond}). We also give a characterization of bounded co-$t$-structures in terms of relative homological algebra. 
The relationship between silting classes and co-$t$-structures is also studied. We prove that a silting class $\omega$ induces a bounded non-degenerated co-$t$-structure on the smallest thick triangulated subcategory of $\T$ containing $\omega.$ We also give a description of the bounded co-$t$-structures on $\T$ (see Theorem \ref{Msc}). Finally, as an application to the particular case of the bounded derived category $\D(\HH),$ where $\HH$ is an abelian
hereditary category which is Hom-finite, Ext-finite and has a tilting object (see \cite{HR}), we give a bijective correspondence
between finite silting generator sets $\omega=\add\,(\omega)$ and bounded co-$t$-structures (see Theorem \ref{teoH}).
\end{abstract}
\maketitle
\section{Introduction.}

In \cite{H}, Hashimoto defined the ``Auslander-Buchweitz context" for abelian categories, giving a new framework
to homological approximation theory. The starting point of Hashimoto's work is the theory of approximations in abelian categories developed by Auslander and Buchweitz in \cite{AB}, which has been a starting point for performing relative homological algebra with respect to suitable subcategories, with applications ranging from the study of Cohen-Macaulay modules over commutative 
rings, to tilting theory, the theory of quasi-hereditary algebras and reductive groups, the study 
of homological conjectures for finite dimensional algebras, and many other topics. On the other hand, in \cite{Be},  Beligiannis generalizes to exact categories the fundamental work of \cite{AB}. In particular, following Hashimoto's ideas, he introduces the Auslander-Buchweitz context for exact categories, which are more general than abelian ones.

In the case of $\mathrm{mod}\,(\Lambda)$ (the category of finitely generated modules over
an artin algebra $\Lambda$), it is important to mention the work of
Auslander and Reiten in \cite{AR}. They studied the notion of approximations of modules using tilting
and cotilting modules, and showed that there is a bijective correspondence
between the basic cotilting modules in $\mathrm{mod}\,(\Lambda)$,
and certain precovering subcategories $\mathcal{X}$ of $\mathrm{mod}\,(\Lambda).$ The main aim in \cite{AR} is to explore the connection between various aspects of tilting theory and the theory of cotorsion pairs in $\mathrm{mod}\,(\Lambda).$

As we mentioned before in \cite{MSSS1}, abelian categories used to be the proper context for the study of homological algebra. But recently, triangulated categories entered into the subject in a relevant way. In \cite{MSSS1}, an analogue of Auslander-Buchweitz approximation theory is developed.

The main aim of the present paper is to explore, in the setting described in \cite{MSSS1}, 
results analogous to the results of Auslander-Reiten in connection with various aspects of tilting theory and the theory of co-$t$-structures. To do that, we use the notions and machinery of \cite{MSSS1}, concentrating our study to the
relations between Auslander-Buchweitz contexts in a triangulated category $\T$ and co-$t$-structures defined on $\T.$

The notion of co-$t$-structure was recently introduced
independently by Pauksztello \cite{P} and Bondarko \cite{Bo} (under the name ``weight structures''). This notion seems to be important, and one reason for this is that they provide important information in a triangulated category
$\T$ allowing the existence of nice ``weight'' decompositions and filtrations. Furthermore, co-$t$-structures provide examples of torsion theories in Krull-Schmidt triangulated categories in the sense of Iyama and Yoshino \cite{Iyama}.

Throughout this paper, $\mathcal{T}$ denotes an arbitrary triangulated category. Given a class $\mathcal{X}$ of objects
of $\mathcal{T},$ the smallest triangulated (respectively, smallest thick) subcategory of $\T$ containing $\X$ is denoted by $\Delta_\T(\X)$ (respectively, $\overline{\Delta}_\T(\X)$).

The paper is organized as follows. In Section 1, we recall, from \cite{MSSS1}, some notions about the Auslander-Buchweitz approximation theory that will be useful in this paper.

In Section 2, we show that the notion of
 relative Auslander-Buchweitz context for triangulated categories $\T$
coincides with the notion of co-$t$-structure on $\overline{\Delta}_\T(\X)$ (see Theorem \ref{M2}). In particular, an Auslander-Buchweitz context is the same as a bounded below co-$t$-structure.
Moreover, we establish a bijective correspondence between the relative
Auslander-Buchweitz contexts $(\mathcal{X},\mathcal{Y})$ on $\mathcal{T}$ and the class of pairs $(\mathcal{X},\omega)$ such that $\mathcal{X}$ is cosuspended and $\omega$ is an $\X$-injective
weak-cogenerator in $\mathcal{X}$ (see Theorem \ref{correspond}).

In Section 3, we focus our attention on bounded, faithful and non-degenerate co-$t$-structures. A characterization of bounded co-$t$-structures, in terms of relative homological algebra, is also given. Furthermore, a relationship between the different types of co-$t$-structures is also established (see Theorem \ref{equivTipos}). We also provide, on one hand, a relationship between several subcategories attached to co-$t$-structures; and on the other hand, some relations between relative homological dimensions. We finish the section with some results involving co-t-structures and the notion of categorical cogenerator.

In Section 4, we study the relationship between co-$t$-structures and silting
classes. In this section, we establish a bijective correspondence between silting classes in $\T$ and bounded co-$t$-structures on the thick subcategory of $\T$ generated by the silting class
(see Corollary \ref{siltingcorresp}). Furthermore, we give a characterization of the bounded 
co-$t$-structures on $\T$ (see Theorem \ref{Msc}). 

In Section 5, we apply the results, obtained in Section 4, to the particular case of the bounded derived category $\D(\HH)$ where $\HH$ is an abelian
hereditary category which is Hom-finite, Ext-finite and has a tilting object. We give a bijective 
correspondence between finite silting generator sets $\omega=\add\,(\omega)$ and 
bounded co-$t$-structures (see Theorem \ref{teoH}). As a nice consequence, we get that any 
bounded co-$t$-structure on $\D(\HH)$ has two companions as $t$-structures: one on the left and 
the other on the right. That is, any bounded co-$t$-structure on $\D(\HH)$ is always left 
(respectively, right) adjacent to a $t$-structure on $\D(\HH)$ in the sense of \cite{Bo}

Note that in \cite{Bo}, the author studies co-$t$-structures on
triangulated categories with arbitrary coproducts (his notion of ``negative
subcategories'' correspond to our notion of silting). In this context, he proves that any silting subcategory $\omega$ provides a co-$t$-structure on the smallest triangulated subcategory of $\T$ closed under arbitrary coproducts and containing  $\omega.$ Our result (Theorem \ref{siltingteo}), which is proved using relative homology techniques, is the analogue for thick subcategories
containing $\omega$, to the Theorem 4.3.2 in \cite{Bo} which was
proved with different techniques.

\section{Preliminaries}
\

Throughout this paper, $\T$ will be a triangulated category and $[1]:\T\rightarrow\T$ its
 suspension functor. 
\

In this paper, when we say that $\mathcal{C}$ is a subcategory of $\T,$ it always means that $\mathcal{C}$ is  a full subcategory which is additive and closed under isomorphisms. For a class $\X$ of objects of $\T,$ we denote by $\add\,(\X)$ the smallest subcategory of $\T$ containing $\X,$ closed under finite direct sums and direct summands.
\

For some classes $\X$ and $\Y$ of objects in $\T,$ we write
${}^\perp\X:=\{Z\in\T\,:\,\Hom_\T(Z,-)|_{\X}=0\}$ and $\X^\perp:=\{Z\in\T\,:\,\Hom_\T(-,Z)|_{\X}=0\}.$
\noindent We also recall that $\X*\Y$ denotes the class of objects $Z\in\T$ for which there exists a distinguished triangle $X\rightarrow Z\rightarrow Y\rightarrow X[1]$ in $\T$ with $X\in\X$ and $Y\in\Y.$ Furthermore, it is said that  $\X$ is {\bf closed under extensions} if $ \X*\X\subseteq \X.$
\

Recall that  a class $\X$ of objects in $\T$ is said to be {\bf suspended} (respectively,  {\bf cosuspended}) if $\X[1]\subseteq\X$ (respectively,  $\X[-1]\subseteq\X$) and $\X$ is  closed under
extensions. Observe that a suspended (respectively, cosuspended) class $\X,$ of objects in $\T,$ is a  subcategory of $\T$ (see 
\cite[Lemma 2.1 (b)]{MSSS1}).
\

Given a class $\X$ of objects in $\T,$ it is said that $\X$ is {\bf{closed under cones}} if for any distinguished triangle  $A\rightarrow B\rightarrow C\rightarrow A[1]$ in $\T$ with  $A,B\in \X$ we have that $C\in \X.$  Similarly, $\X$ is {\bf{closed under cocones}} if for any distinguished triangle  $A\rightarrow B\rightarrow C\rightarrow A[1]$ in $\T$ with $B,C\in \X$ we have that $A\in \X.$
\

Let $\mathcal{X}$ be a class of objects of $\mathcal{T}.$ We denote by $\U_{\X}$ (respectively,  ${}_{\X}\U$) the smallest suspended (respectively,  cosuspended) subcategory of $\T$ containing the class $\X.$ Note that if $\X$ is suspended (respectively, cosuspended) subcategory of $\T,$ then $\X=\U_{\X}$ (respectively,  $\X={}_{\X}\U$). We also recall that a subcategory $\U$ of $\T,$ which is suspended and cosuspended, is called a {\bf{triangulated subcategory}} of $\T.$ A {\bf{thick}} subcategory of $\T$ is a triangulated subcategory of $\T$ which is closed under direct summands in $\T.$  The smallest triangulated (respectively, smallest thick) subcategory of $\T$ containing $\X$ is denoted by $\Delta_\T(\X)$ (respectively, $\overline{\Delta}_\T(\X)$).
\

We recall the following well known definition (see, for example, \cite{B} and \cite{BR}).

\begin{defi} Let $\X$ and $\Y$ be classes of objects in a triangulated category $\T.$ A morphism $f:X\to C$ in $\T$ is said to be an {\bf $\X$-precover} of $C$ if $X\in\X$ and $\Hom_\T(X',f):\Hom_\T(X',X)\to\Hom_\T(X',C)$ is surjective $\forall X'\in\X.$ If any $C\in\Y$ admits an $\X$-precover, then $\X$ is called a {\bf precovering} class in $\Y.$ By dualizing the definition above, we get the notion of an {\bf $\X$-preenveloping} of $C$ and a {\bf preenveloping} class in $\Y.$ Finally, it is said that $\X$ is {\bf{functorially finite}} in $\T$ if $\X$ is both precovering and preenveloping in $\T.$
\end{defi}

Now, we recall from \cite{MSSS1}, the following definitions. For a more completed discussion and properties of such notions, we suggest that the reader see \cite{MSSS1}.

\begin{defi}\cite{MSSS1}\label{epsilon} Let $\X$ be a class of objects in $\T.$ For any natural number $n,$ we introduce inductively the  
class $\varepsilon^\wedge_n(\X)$ as follows: $\varepsilon^\wedge_0(\X):=\X$ and assuming defined $\varepsilon^\wedge_{n-1}(\X),$ the class 
$\varepsilon^\wedge_{n}(\X)$ is given by all the objects $Z\in\T$ for which there exists a distinguished triangle in $\T$
$$\begin{CD}
Z[-1] @>>> W @>>> X @>>> Z \,
\end{CD}$$ with $W\in\varepsilon^\wedge_{n-1}(\X)$ and $X\in\X.$

Dually, we set
 $\varepsilon^\vee_0(\X):=\X$ and assuming defined $\varepsilon^\vee_{n-1}(\X),$ the class $\varepsilon^\vee_{n}(\X)$ is formed for all the 
objects $Z\in\T$ for which there exists a distinguished triangle in $\T$
$$\begin{CD}
Z @>>> X @>>> K @>>>  Z[1]\,
\end{CD}$$ with $K\in\varepsilon^\vee_{n-1}(\X)$ and $X\in\X.$ We also introduce the 
following classes
$$\X^\wedge:=\cup_{n\geq 0}\;\varepsilon^\wedge_n(\X),\quad\X^\vee:=\cup_{n\geq 0}\;\varepsilon^\vee_n(\X)\text{ and }\X^\sim :=(\X^\wedge)^\vee.$$
\end{defi}

For the convenience of the reader, we include the following remark from \cite{MSSS1}.

\begin{rk}\cite[Remark 3.6 (2)]{MSSS1}  \label{varios}
Let $(\Y,\omega)$ be a pair of classes of objects in $\T$ with $\omega\subseteq\Y.$ If $\Y$
is closed under cones (respectively, cocones) then $\omega^\wedge\subseteq\Y$ (respectively,
 $\omega^\vee\subseteq\Y$). Indeed, assume that $\Y$ is closed under cones and let
$M\in\omega^\wedge.$ Thus $M\in\varepsilon^\wedge_n(\omega)$ for some $n\in\Bbb{N}.$ If $n=0$
 then $M\in\omega\subseteq \Y.$ Let $n>0,$ and hence there is a distinguished triangle 
$M[-1] \to K\to Y\to M$ in $\T$ with $K\in\varepsilon^\wedge_{n-1}(\omega)$ and $Y\in\Y.$ By 
induction $K\in\Y$ and hence $M\in\Y$ since $\Y$ is closed under cones; proving that 
$\omega^\wedge\subseteq\Y.$
\end{rk}

In what follows, to deal with the (co) resolution, relative projective and relative injective 
dimensions, we consider the extended natural numbers 
$\overline{\mathbb{N}}:=\mathbb{N}\cup\{\infty\}.$ Here, we set the following rules: 
(a) $x+\infty=\infty$ for any $x\in\overline{\mathbb{N}},$ and (b) $x<\infty$ for any 
$x\in\mathbb{N}.$ Finally, we declare, by definition, that the minimum of the empty set is 
 $\infty.$ That is, $\mini(\emptyset):=\infty.$

\begin{defi}\cite{MSSS1} Let $\X$ be a class of objects in $\T,$ and let $M\in\T.$
The $\X$-{\bf{resolution dimension}} of $M$ is $$\resdim_{\X}(M):=\mini\,\{n\in\mathbb{N} 
\;:\; M\in\varepsilon^\wedge_{n}(\X) \}.$$ Dually, the $\X$-{\bf{coresolution dimension}} of 
$M$ is $$\coresdim_{\X}(M):=\mini\,\{n\in\mathbb{N} \;:\; M\in\varepsilon^\vee_{n}(\X) \}.$$
\end{defi}

The following result, and its dual version, will be used in Section 2.

\begin{teo}\cite[Theorem 3.5]{MSSS1}\label{M3(a)} 
 For any cosuspended subcategory $\X$ of $\T$ and any object $C\in\T,$  the following statements hold.
 \begin{itemize}
  \item[(a)] $\resdim_\X(C)\leq n$ if and only if $C\in\X[n].$
  \item[(b)] $\X^\wedge=\cup_{n\geq 0}\,\X[n]=\Delta_\T(\X).$
  \item[(c)] If $\X$ is closed under direct summands in $\T,$ then $\X^\wedge=\overline{\Delta}_\T(\X).$
 \end{itemize}
\end{teo}

For the convenience of the reader, we include the dual version of \ref{M3(a)}.

\begin{rk}\label{M3(a)dual} For any suspended subcategory $\Y$ of $\T$ and any object $C\in\T,$  the following statements hold.
 \begin{itemize}
  \item[(a)] $\coresdim_\Y(C)\leq n$ if and only if $C\in\Y[-n].$
  \item[(b)] $\Y^\vee=\cup_{n\geq 0}\,\Y[-n]=\Delta_\T(\Y).$
  \item[(c)] If $\Y$ is closed under direct summands in $\T,$ then $\Y^\vee=\overline{\Delta}_\T(\Y).$
 \end{itemize}
\end{rk}

We recall the notion of $\X$-projective (respectively, $\X$-injective) dimension of objects in $\T.$

\begin{defi}\cite{MSSS1} Let $\X$ be a class of objects in $\T$ and $M$ an object in $\T.$
\begin{enumerate}
 \item[(a)] The {\bf{$\X$-projective dimension}} of
$M$ is
$$\pd_{\X}(M):=\mini\,\{n\in\mathbb{N}\; :\; \Hom_\T(M[-i],-)\mid_{\X}=0, \quad\forall
i>n\}.$$
 \item[(b)] The {\bf{$\X$-injective dimension}} of $M$ is
$$\id_{\X}(M):=\mini\,\{n\in\mathbb{N}\; :\; \Hom_\T(-,M[i])\mid_{\X}=0, \quad\forall i>n\}.$$
 \end{enumerate}
\end{defi}

\begin{defi}\cite{MSSS1} Let $(\X,\omega)$ be a pair of classes of objects in $\T.$
\begin{itemize}
 \item[(a)] $\omega$ is a {\bf weak-cogenerator} in $\X,$ if $\omega\subseteq\X\subseteq\X[-1]*\omega.$
 \item[(b)] $\omega$ is a {\bf weak-generator} in $\X,$ if $\omega\subseteq\X\subseteq \omega * \X[1].$
 \item[(c)] $\omega$ is $\X$-{\bf{injective}} if $\id_{\X}(\omega)=0;$ and dually, $\omega$ is $\X$-{\bf{projective}} if $\pd_{\X}(\omega)=0.$
\end{itemize}
\end{defi}

\section{Relative Auslander-Buchweitz context and co-$t$-structures}

In this section, we give the notion of the (relative) Auslander-Buchweitz context for a triangulated category $\T,$ relating this notion with the concept of co-$t$-structure.

\begin{defi}\label{defcot}\cite{Bo, P} A pair $(\A,\B)$ of subcategories in $\T$ is said to be a {\bf{co-$t$-structure}} on $\T$ if the following conditions hold.
 \begin{itemize}
  \item[(a)] $\A$ and $\B$ are closed under direct summands in $\T.$
  \item[(b)] $\A[-1]\subseteq\A$ and $\B[1]\subseteq \B.$
  \item[(c)] $\Hom_\T(\A[-1],\B)=0.$
  \item[(d)] $\T=\A[-1]*\B.$
 \end{itemize}
\end{defi}

We will make use of the following result, stated by D. Pauksztello in \cite{P}.
\begin{pro}\label{proP}\cite[Proposition 2.1]{P} Let $(\A,\B)$ be a co-$t$-structure on $\T.$ Then, the following statements hold.
 \begin{itemize}
 \item[(a)] $\A[-1]$ is a precovering class in $\T.$
 \item[(b)] $\B$ is a preenveloping class in $\T.$
 \item[(c)] $\A[-1]={}^\perp\B$ and $\B=\A^\perp[-1].$
 \item[(d)] $\A$ and $\B$ are closed under extensions.
 \end{itemize}
\end{pro}

\begin{lem}\label{idsubseteq} Let $(\A,\B)$ be a co-$t$-structure on $\T,$ and
$\Y$ be a class of objects in $\T.$ Then, the following statements hold.
\begin{itemize}
 \item[(a)] $\id_\A(\Y)\leq n$ if and only if $\Y\subseteq\B[-n].$
 \item[(b)] $\pd_\B(\Y)\leq n$ if and only if $\Y\subseteq\A[n].$
\end{itemize}
\end{lem}
\begin{dem} (a) By \cite[Lemma 4.2]{MSSS1}, we get the equivalence: $\id_\A(\Y)\leq n$ if and 
only if $\Y\subseteq{}_\A\U^\perp[-n-1].$ Therefore, since $(\A, \B)$ is a co-$t$-structure, it follows from \ref{proP} (c), that ${}_\A\U^\perp[-n-1]=\B[-n].$
\

(b) It is dual to (a).
\end{dem}
${}$
\vspace{.2cm}

The following result states that, for a co-$t$-structure $(\A,\B)$ on $\T,$ the class $\omega:=\A\cap\B$ is an $\A$-injective weak-cogenerator in $\A;$ and moreover, $\omega$ is also a $\B$-projective weak-generator in $\B.$ Note that $\omega=\add\,(\omega).$

\begin{pro}\label{cogenXiy} Let $(\A,\B)$ be a co-$t$-structure on $\T,$ and let $\omega:=\A\cap\B.$ Then, the following statements hold.
 \begin{itemize}
  \item[(a)] $\id_\A(\B)=0\quad\text{and}\quad \A\subseteq\A[-1]*\omega.$
  \item[(b)] $\pd_\B(\A)=0\quad\text{and}\quad \B\subseteq\omega*\B[1].$
 \end{itemize}
\end{pro}
\begin{dem} By \ref{idsubseteq}, it follows that $\id_\A(\B)=0=\pd_\B(\A).$ To see the inclusion 
in (a), let $C\in\A.$ Then, by \ref{defcot} (d), we have a distinguished triangle $C'\to C\to C''\to C'[1]$ in $\T$ with $C'\in\A[-1]$ and $C''\in\B.$ Hence, by \ref{proP} (d), it follows that $C''\in\A\cap\B=\omega;$ proving that $\A\subseteq\A[-1]*\omega.$
\

By \ref{defcot}, for any $X\in\T,$ there is a distinguished triangle $A\to X\to B[1]\to A[1]$ in $\T$ with $A\in\A$ and $B\in\B.$ Moreover, in the case $X\in\B,$ it follows from the preceding triangle that $A\in\A\cap\B=\omega;$ getting us that $\B\subseteq\omega*\B[1].$
\end{dem}
${}$
\vspace{.2cm}

We will show the relation between the notions of  cosuspended (respectively, suspended) 
subcategories $\X$, weak-cogenerator (respectively, weak-generator), $\X$-injective 
(respectively,  $\X$-projective) and co-$t$-structures on $\overline{\Delta}_\T(\X).$ We only 
state the results for the cosuspended case and omit those for the suspended case which can be proved by similar arguments.
\

First, we show that  any  $\X$-injective weak-cogenerator  in a cosuspended subcategory 
$\X=\add\,(\X)$ of $\T$ provides a co-$t$-structure on $\overline{\Delta}_\T(\X)=\X^\wedge.$

\begin{teo}\label{M1} Let $(\X,\omega)$ be a pair of classes of objects in $\T$ which are closed under direct summands, $\X$ be cosuspended and $\omega$ be an $\X$-injective weak-cogenerator in $\X.$ Then, the following statements hold.
\begin{enumerate}
 \item[(a)] The pair $(\X^\wedge\cap{}^\perp(\omega^\wedge)[1],\omega^\wedge)$
is a co-t-structure on the triangulated category $\X^\wedge.$
 \item[(b)] $\omega^\wedge=\X^\wedge\cap\X^\perp[-1],\quad$ $\X=\X^\wedge\cap{}^\perp(\omega^\wedge)[1]\quad$ and $\quad\omega=\X\cap\X^\perp[-1].$
 \item[(c)] If $\omega'$ is an $\X$-injective weak-cogenerator in $\X,$ then $\omega=\add\,\omega'.$
\end{enumerate}
\end{teo}
\begin{dem} First note that $ {\X}= {}_\X\U$ since $\X$ is cosuspended.
\

(b) From  \cite[Proposition 5.9]{MSSS1}, we have that $\omega^\wedge=\X^\wedge\cap\X^\perp[-1]$ . By \cite[Proposition 5.2 (b)]{MSSS1}, it 
follows that $\omega=\X\cap\X^\perp[-1]$ since $\X$ is cosuspended. Moreover, by \cite[Theorem 5.10]{MSSS1} it follows that 
$\X=\X^\wedge\cap{}^\perp(\omega^\wedge)[1].$
\

(a)  We have that $\omega^\wedge=\X^\wedge\cap\X^\perp[-1]$ is suspended and closed under direct summands. Therefore $\X^\wedge\cap{}^\perp(\omega^\wedge)[1]$ is cosuspended and closed under direct summands. So, in order to get that the given pair in (a) is a co-$t$-structure on the triangulated category $\X^\wedge,$ it is enough to see that $\X^\wedge=(\X^\wedge\cap{}^\perp(\omega^\wedge))*\omega^\wedge.$ But this is a consequence of \cite[Corollary 5.5 (b)]{MSSS1} since $\X[-1]=\X^\wedge\cap{}^\perp(\omega^\wedge).$
\

(c) It follows from (b) and the fact that $\add\,(\omega')$ is an $\X$-injective weak-cogenerator in $\X.$
\end{dem}

\begin{rk} \label{wcg}Let  $\X=\add\,(\X)$ be a cosuspended subcategory of $\T.$   Note that $\X\cap\X^\perp[-1]$ is $\X$-injective. Moreover, from \ref{M1}, we get that: If there is an $\X$-injective weak-cogenerator  $\omega=\add\,(\omega)$ in $\X$ then it is unique.  Consequently,
   there is an $\X$-injective weak-cogenerator  $\omega=\add\,(\omega)$ in $\X$  if and only if  $\X\cap\X^\perp[-1]$ is a weak-cogenerator in $\X.$
\end{rk}

The Auslander-Buchweitz context for abelian categories was introduced by M. Hashimoto in 
\cite{H}. Inspired by that, we will introduce such a context for a triangulated category 
$\T.$ To do so, we define the notion of a relative Auslander-Buchweitz context on $\T.$ 
Observe that the ``relative Auslander-Buchweitz context'' in triangulated categories is used 
 for an analogue of what Hashimoto calls ``weak Auslander-Buchweitz context'' in abelian 
 categories.

\begin{defi} Let $(\X,\Y)$ be a pair of classes of objects in $\T,$ and let $\omega:=\X\cap\Y.$ The pair $(\X,\Y)$ is said to be a {\bf{relative Auslander-Buchweitz context}} on $\T$ if the following three conditions hold:
   \begin{itemize}
    \item[(AB1)] $\X$ is cosuspended and closed under direct summands in $\T.$
    \item[(AB2)] $\Y$ is suspended and closed under direct summands in $\T$ and $\Y\subseteq\X^\wedge.$
    \item[(AB3)] $\omega$ is an $\X$-injective weak-cogenerator in $\X.$
   \end{itemize}
The pair $(\X,\Y)$ is said to be an {\bf{Auslander-Buchweitz context}} on $\T$ if $(\X,\Y)$ is a relative Auslander-Buchweitz context on $\T$ and $\X^\wedge=\T.$
\end{defi}

\begin{teo}\label{M2} Let $(\X,\Y)$ be a relative Auslander-Buchweitz context on $\T$ and $\omega:=\X\cap\Y.$ Then, the following statements hold.
 \begin{itemize}
 \item[(a)] $\omega=\X\cap\X^\perp[-1]\quad{and}\quad\omega^\wedge=\Y.$
 \item[(b)] $(\X,\Y)$ is a co-$t$-structure on the triangulated category $\X^\wedge.$
\end{itemize}
\end{teo}
\begin{dem} (a) The first equality follows from \ref{M1}. Since $\omega\subseteq\Y$ and $\Y$ is suspended, it follows from \ref{varios} that $\omega^\wedge\subseteq\Y.$\\
We assert that $\id_\X(\Y)=0.$ Indeed, let $C\in\Y\subseteq\X^\wedge.$ Hence, by \cite[Theorem 5.4]{MSSS1}, we have a distinguished triangle $Y_C\to X_C\to C\to Y_C[1]$ in $\T$ with $X_C\in\X$ and $Y_C\in\omega^\wedge\subseteq\Y.$ Hence $X_C\in\X\cap\Y=\omega$ and so $\id_\X(X_C)=0.$ On the other hand, since $\id_\X(Y_C)=0$ (see \cite[Proposition 5.2 (a)]{MSSS1}), it follows by \cite[Lemma 5.7]{MSSS1} that $\id_\X(C)=0;$ proving the assertion. Finally, $\id_\X(\Y)=0$ and the fact that $\X$ is cosuspended implies by \cite[Lemma 4.2]{MSSS1} that $\Y\subseteq\X^\wedge\cap\X^\perp[-1].$ Therefore $\Y\subseteq\omega^\wedge$ by \ref{M1}.
\

(b) Since $\omega^\wedge=\Y,$ we have that (b) follows from \ref{M1}.
\end{dem}
${}$
\vspace{.2cm}

Given a class $\X$ of objects in $\T,$ we recall that $\overline{\Delta}_\T(\X)$ denotes the smallest thick subcategory of $\T$ containing the class $\X.$

\begin{pro}\label{M3} Let $\X$ and $\Y$ be classes of objects in $\T$ such that the pair $(\X,\Y)$ is a co-$t$-structure on the triangulated category $\overline{\Delta}_\T(\X).$ Then, $\overline{\Delta}_\T(\X)=\X^\wedge$ and $(\X,\Y)$ is a relative Auslander-Buchweitz context on $\T.$
\end{pro}
\begin{dem} By \ref{proP} (d), we have that $\X$ is cosuspended and $\Y$ is suspended. In particular, from \ref{M3(a)}, we conclude that 
$\overline{\Delta}_\T(\X)=\X^\wedge.$ The fact that $\omega=\X\cap\Y$ is an $\X$-injective weak-cogenerator in $\X,$ follows from 
\ref{cogenXiy} (a).
\end{dem}

Now, we are in a position to state our main result  in this section. In order to do that,  we introduce the following classes.

\begin{defi} For a given triangulated category $\T,$ we introduce the following classes:
\begin{itemize}
 \item[(a)] $\C_1$ consists of all pairs $(\X,\omega)$ of classes of objects in $\T,$ which are closed under direct summands, and such that $\X$ is cosuspended and $\omega$ is an $\X$-injective weak-cogenerator in $\X.$
 \item[(b)] $\C_2$ consists of all pairs $(\X,\Y)$ of classes of objects in $\T,$ which are a relative Auslander-Buchweitz context on $\T.$
 \item[(c)] $\C_3$ consists of all pairs $(\X,\Y)$ of classes of objects in $\T,$ which are a co-$t$-structure on
$\overline{\Delta}_\T(\X).$
 \item[(d)] $\C_4$ consists of all cosuspended subcategories $\X$ in $\T,$ which are  precovering in $\X^\wedge$ and $\X=\add\,(\X).$
\end{itemize}
\end{defi}

An additive category $\C$ is said to be Krull-Schmidt if any object $C$ in $\C$ has a finite decomposition $C=C_1\oplus C_2\oplus\cdots\oplus C_n$ such that each $C_i$ is indecomposable with local endomorphism ring.

Let $R$ be a commutative artinian ring. We recall that an $R$-linear triangulated category $\T$ is said to be $\Hom$-finite if $\Hom_\T(X,Y)$ if a finite generated $R$-module for any $X,Y\in\T.$

Let $\Lambda$ be an artin $R$-algebra. It is well known that the bounded derived category $\Kb$ is a typical example of an $R$-linear triangulated category which is Krull-Schmidt and $\Hom$-finite.

\begin{teo}\label{correspond} Let $\T$ be a triangulated category. Then, the following statements hold.
\begin{itemize}
 \item[(a)] $\C_2=\C_3$ and the correspondence $\C_1\to \C_2,$ $(\X,\omega)\mapsto (\X,\Y:=\omega^\wedge),$ is a bijection with inverse $\C_2\to \C_1$ given by $(\X,\Y)\mapsto (\X,\omega:=\X\cap\Y).$
 \item[(b)] If $\T$ is an $R$-linear triangulated category which is $\Hom$-finite and Krull-Schmidt, then  the correspondence $\C_4\to \C_3,$
$\X\mapsto (\X,\Y:=\X^\perp[-1]\cap\X^\wedge)$ is a bijection with inverse $\C_4\to \C_3$ given by $(\X,\Y)\mapsto \X.$
\end{itemize}
\end{teo}
\begin{dem} (a) It follows from \ref{M1}, \ref{M2} and \ref{M3}.
\

(b) Assume that $\T$ is a $R$-linear triangulated category which is $\Hom$-finite and Krull-Schmidt. Let $\X\in\C_4.$ Since $\X$ is cosuspended and closed under 
direct summands in $\T,$ it follows from \ref{M3(a)} that $\overline{\Delta}_\T(\X)=\X^\wedge.$ On the other hand, by 
\cite[Proposition 2.3]{Iyama}, we have that $\Hom_\T(\X,\X^\perp\cap\X^\wedge)=0$  and $\X^\wedge=\X*(\X^\perp\cap\X^\wedge).$ Therefore 
$$\X[-1]*(\X^\perp[-1]\cap\X^\wedge)=(\X*(\X^\perp\cap\X^\wedge))[-1]=\X^\wedge[-1]=\X^\wedge,$$ getting us that 
$(\X,\X^\perp[-1]\cap\X^\wedge)\in\C_3.$
\

Consider a pair $(\X,\Y)\in\C_3.$ Then by \ref{M3}, \ref{M2} and \ref{M1}, it follows that $\Y=\X^\perp[-1]\cap\X^\wedge.$ Moreover, since the pair $(\X[1],\Y[1])$ is also a co-$t$-structure on $\overline{\Delta}_\T(\X),$ we have from \ref{proP} that $\X\in\C_4.$ Furthermore, since $\Y=\X^\perp[-1]\cap\X^\wedge,$ it follows that the correspondence $\X\mapsto (\X,\X^\perp[-1]\cap\X^\wedge)$ induces a bijection, with inverse $(\X,\Y)\mapsto \X,$ between the classes $\C_4$ and $\C_3.$
\end{dem}

\begin{cor} \label{wcgext} There is a bijective correspondence $\X\mapsto (\X,\overline{\Delta}_\T(\X)\cap\X^\perp[-1])$ between cosuspended subcategories $\X=\add\,(\X)$ of $\T$ such that $\X\cap\X^\perp[-1]$ is a weak-cogenerator in $\X,$  and co-$t$-structures  $(\X,\Y)$ on $\overline{\Delta}_\T(\X).$
  \end{cor}
\begin{dem} It follows from \ref{correspond} and   \ref{wcg}.
\end{dem}

\section{Bounded, faithful and non-degenerate co-$t$-structures}

In this section we focus our attention on bounded, faithful and non-degenerate co-t-structures. We finish the section with some results involving co-t-structures and the notion of categorical cogenerator.
\

Following the terminology for co-$t$-structures on triangulated categories given in \cite{Bo}, we recall the following definition.

\begin{defi} Let $(\A, \B)$ be a co-$t$-structure on  $\T.$ It is said that
 $(\A, \B)$ is {\bf{ bounded below }} (respectively, {\bf{ bounded above}})  if $\cup_{n\in\Bbb{Z}}\,\A[n]=\T$ (respectively, 
$\cup_{n\in\Bbb{Z}}\,\B[n]=\T$). So, the pair  $(\A, \B)$ is said to be {\bf{ bounded  }} if it is bounded both below and above.
\end{defi}

\begin{rk}\label{bog} From \ref{M3(a)} and \ref{M3(a)dual}, we have that a co-$t$-structure $(\A, \B)$ on  $\T$ is  bounded below 
(respectively, above)  if and only if $\A^\wedge=\T$ (respectively, $\B^\vee=\T$).
\end{rk}

\begin{cor} \label{bwcgext} There is a bijective correspondence $\X\mapsto (\X,\X^\perp[-1])$ between cosuspended subcategories 
$\X=\add\,(\X)$ of $\T$ such that $\X^\wedge=\T$ and $\X\cap\X^\perp[-1]$ is a weak-cogenerator in $\X,$  and bounded below 
co-$t$-structures  $(\X,\Y)$ on $\T.$
 \end{cor}
\begin{dem} It follows from \ref{wcgext} and \ref{bog}.
\end{dem}
${}$
\vspace{.2cm}

Now, we prove some relationships between the relative homological dimensions attached to a co-$t$-structure.

\begin{pro}\label{relpdresct} Let $(\A,\B)$ be a co-$t$-structure on $\T$ and $\omega:=\A\cap\B.$ Then
\begin{itemize}
 \item[(a)] $\pd_{\B}(M)=\resdim_{\A}(M)$ and $\id_{\A}(M)=\coresdim_{\B}(M),\quad \forall M\in\T.$
 \item[(b)] $\resdim_{\A}(M)=\resdim_\omega(M),\quad \forall M\in\omega^\wedge.$
 \item[(c)] $\coresdim_{\B}(M)=\coresdim_\omega(M),\quad \forall M\in\omega^\vee.$
\end{itemize}
\end{pro}
\begin{dem} By \ref{proP}, we know that $\B=\A^\perp[-1]={}_{\A}\,\U^\perp[-1]$ and $\A={}^\perp\B[1]={}^\perp\U_{\,\B}[1].$ Hence, from \cite[Proposition 4.3]{MSSS1}, we get (a). Finally, (b) and (c) follows from \cite[Theorem 4.4]{MSSS1} and its dual, and the item (a).
\end{dem}
${}$
\vspace{.2cm}

The next result provides a relationship between several subcategories attached to co-$t$-structures. Furthermore, it
characterizes the bounded below co-$t$-structures on $\T.$ We recall that $\omega^\sim:=(\omega^\wedge)^\vee$ for any class $\omega$
of objects in $\T.$

\begin{teo}\label{criterio} Let $(\A, \B)$ be a co-$t$-structure on  $\T$ and $\omega:=\A\cap\B.$ Then, the following conditions hold.
 \begin{itemize}
 \item[(a)] $\U_\omega=\omega^\wedge=\A^\wedge\cap\B\;$ and $\;{}_\omega\,\U=\omega^\vee=\B^\vee\cap\A.$
 \item[(b)] $\overline{\Delta}_\T(\omega)=\omega^\sim=\{C\in\A^\wedge : \id_{\A}(C)<\infty\}=\{C\in\B^\vee : 
\pd_{\B}(C)<\infty\}=\A^\wedge\cap\B^\vee.$
 \item[(c)] The following conditions are equivalent:
   \begin{itemize}
    \item[(c1)]  $(\A, \B)$  is bounded below.
    \item[(c2)]  $\B\subseteq\omega^\sim.$
    \item[(c3)]  $\omega^\wedge=\B.$
    \item[(c4)]  $\B\subseteq\A^\wedge.$
   \end{itemize}
 \end{itemize}
\end{teo}
\begin{dem} (a) Since  $(\A, \B)$ is a co-$t$-structure on $\T,$ we obtain from \ref{cogenXiy} that $\omega$ is an $\A$-injective 
weak-cogenerator in $\A.$ Therefore, the first equalities in (a) follows from \cite[Proposition 5.9]{MSSS1}, and the second ones can be 
proven by dualizing \cite[Proposition 5.9]{MSSS1}.
\

(b) It follows from \cite[Theorem 5.16]{MSSS1} and its dual. 
\

(c) (c1) $\Rightarrow$ (c3) Let $\A^\wedge=\T$ (see \ref{bog}). Then, by \ref{M3}, it follows that $(\A, \B)$ is an Auslander-Buchweitz 
context on $\T.$ Hence $\B=\omega^\wedge$ by \ref{M2}.
\

(c3) $\Rightarrow$ (c2) Assume that  $\B=\omega^\wedge.$ Since $\omega^\wedge\subseteq\omega^\sim,$ we get that $\B\subseteq\omega^\sim.$
\

(c2) $\Rightarrow$ (c1) Suppose that $\B\subseteq\omega^\sim.$ We assert that $\T=\A^\wedge.$ Indeed, since $(\A, \B)$ is a 
co-$t$-structure on $\T,$ we have that $\T=\A[-1]*\A^\perp[-1]=(\A*\A^\perp)[-1];$ and so $\T=\A*\A^\perp.$ Thus for any 
$C\in\T$ there is a distinguished triangle $Z[-1]\to A\to C\to Z$ in $\T$ with $A\in\A$ and $Z\in\A^\perp.$ But 
$Z[-1]\in\A^\perp[-1]=\B\subseteq\omega^\sim\subseteq\A^\wedge$ by (b); proving that $C\in\A^\wedge.$
\

(c3) $\Leftrightarrow$ (c4) It follows from the equality $\omega^\wedge=\A^\wedge\cap\B$ (see (a)).
\end{dem}
${}$
\vspace{.2cm}

The  results for bounded above co-$t$-structures can be stated and proved. To give an example, we give the following 
characterization of bounded above co-$t$-structures.

\begin{rk}\label{criterioba} Let $(\A, \B)$ be a co-$t$-structure on  $\T$ and $\omega:=\A\cap\B.$ Then, the following conditions are equivalent:
   \begin{itemize}
    \item[(a)]  $(\A, \B)$  is bounded above.
    \item[(b)]  $\A\subseteq\omega^\sim.$
    \item[(c)]  $\omega^\vee=\A.$
    \item[(d)]  $\A\subseteq\B^\vee.$
   \end{itemize}
\end{rk}

Following the terminology for $t$-structures on triangulated categories, and also \cite{Bo} and \cite{P}, we give the following definitions.

\begin{defi} Let $(\A, \B)$ be a co-$t$-structure on  $\T,$ and let $\omega:=\A\cap\B.$ It is said that $(\A, \B)$ is {\bf{faithful below}} (respectively, {\bf{faithful above}}) if $\;\cup_{n\in\Bbb{Z}}\,\A[n]=\overline{\Delta}_\T(\omega)$ (respectively, $\;\cup_{n\in\Bbb{Z}}\,\B[n]=\overline{\Delta}_\T(\omega)$). So, it is said that $(\A, \B)$ is {\bf{faithful}} if it is both faithful below and above.
\end{defi}

\begin{pro}\label{criteriobf} Let $(\A, \B)$ be a co-$t$-structure on  $\T,$ and let 
$\omega:=\A\cap\B.$ Then, the following conditions are equivalent. 
\begin{itemize}
 \item[(a)] $(\A, \B)$ is faithful below.
 \item[(b)] $\A^\wedge=\overline{\Delta}_\T(\omega).$
 \item[(c)] $\A^\wedge\subseteq \B^\vee.$
 \item[(d)] $(\A, \B)$ is bounded above.
 \end{itemize}
\end{pro}
\begin{dem} (a) $\Leftrightarrow$ (b) It follows from \ref{M3(a)}. 
\

(b) $\Leftrightarrow$ (c) It follows from \ref{criterio} (b).
\

(c) $\Leftrightarrow$ (d) First, observe that $\A^\wedge\subseteq \B^\vee$ is equivalent to the 
inclusion $\A\subseteq \B^\vee$ since $\B^\vee$ is a triangulated subcategory of $\T$ (see \ref{M3(a)dual}). Therefore, by \ref{criterioba}, we get the result.
\end{dem}

\begin{cor}\label{b=f} Let $(\A, \B)$ be a co-$t$-structure on  $\T.$ Then, $(\A, \B)$ is bounded if and only 
if it is faithful.
\end{cor}
\begin{dem} It follows from \ref{criteriobf} and its dual.
\end{dem}

\begin{teo}\label{corocotbf} Let $(\A, \B)$ be a bounded co-$t$-structure on  $\T,$ and let $\omega:=\A\cap\B.$ Then, the following statements hold.
\begin{itemize}
 \item[(a)] $\overline{\Delta}_\T(\omega)=(\omega^\wedge)^\vee=(\omega^\vee)^\wedge=\T,\;$ $\;\U_\omega=\omega^\wedge=\B\;$ and $\;{}_\omega\U=\omega^\vee=\A.$
 \item[(b)] $\id_\A(C)=\id_\omega(C)=\coresdim_\B(C)<\infty$ for all $C\in\T.$
 \item[(c)] $\pd_\B(C)=\pd_\omega(C)=\resdim_\A(C)<\infty$ for all $C\in\T.$
 \item[(d)] $\coresdim_\B(C)=\coresdim_\omega(C)<\infty$ for all $C\in\A.$
 \item[(e)] $\resdim_\A(C)=\resdim_\omega(C)<\infty$ for all $C\in\B.$
\end{itemize}
\end{teo}
\begin{dem} (a) It follows from \ref{criteriobf}, \ref{criterio} and its dual.
\

(b) Since $\omega^\sim=\T$ (see (a)), we get from \ref{criterio} (b) that $\id_\A(C)<\infty$ for all $C\in\T.$ Thus, (b) follows from \ref{relpdresct} (a) and \cite[Proposition 5.17 (a)]{MSSS1}.
\

(c) Using that $(\A, \B)$ is a co-$t$-structure on  $\A^\wedge=\T,$ we get from \ref{cogenXiy}, that the pair $(\A,\omega)$ satisfies the needed hypothesis in \cite[Theorem 5.6]{MSSS1}; proving (c).
\

(d) and (e) They follow from \ref{relpdresct} since $\omega^\vee=\A$ and $\omega^\wedge=\B.$
\end{dem}
\vspace{.2cm}

Now, we will do one application of \ref{corocotbf} to the so called Rouquier's relative dimension which was introduced in \cite{MSSS1}. Let $\X$ and $\Y$ be classes of objects in a triangulated category $\T.$ Consider the subcategory 
$\left\langle \X\right\rangle:=\add\,(\cup_{i\in\Enteros}\,\X[i])$ and let  
$\X\diamondsuit\Y:=\left\langle\X*\Y\right\rangle.$ Following R. Rouquier in \cite{Ro}, we inductively define $\la \X\ra_0:=0$ and $\la \X\ra_n:=\la \X\ra_{n-1}\diamondsuit\la \X\ra$ for $n\geq 1.$ So, we start with the following definition. 

\begin{defi}\cite[Definition 6.3]{MSSS1} Let $\T$ be a triangulated category, $\X$ a class of objects in $\T$ and $M\in\T.$
The $\X$-dimension of $M$ is 
$$\dim_{\X}(M):=\mini\{n\in\N\text{ such that } M\in\la\X\ra_{n+1}\}.$$
\noindent For a class $\Y$ of objects in $\T,$ we set $\dim_\X(\Y):=\sup\,\{
\dim_\X(Y)\;:\;Y\in\Y\}.$
\end{defi}

\begin{cor} Let $(\A, \B)$ be a bounded co-$t$-structure on  $\T,$ and let $\omega:=\A\cap\B.$ Then
\begin{itemize}
 \item[(a)] $\maxi\,\{\dim_\A(C),\dim_\B(C)\}\leq\dim_\omega(C)$ for all $C\in\T.$
 \item[(b)] $\dim_\A(C)\leq\pd_\omega(C)=\pd_\B(C)<\infty$ for all $C\in\T.$
 \item[(c)] $\dim_\B(C)\leq\id_\omega(C)=\id_\A(C)<\infty$ for all $C\in\T.$
 \item[(d)] $\dim_\omega(C)\leq\id_\omega(C)=\id_\A(C)<\infty$ for all $C\in\A.$
 \item[(e)] $\dim_\omega(C)\leq\pd_\omega(C)=\pd_\B(C)<\infty$ for all $C\in\B.$
\end{itemize}
\end{cor}
\begin{dem} Since $\omega=\A\cap\B,$ (a) follows from \cite[Lemma 6.4 (b)]{MSSS1}. Let $\X$ be any class of objects in $\T$ and $M\in\T.$ Since $\X\subseteq\la\X\ra,$ we get  
by \cite[Proposition 6.6]{MSSS1}, that $\dim_\X(M)\leq\mini\,\{\resdim_{\X}(M),\coresdim_{\X}(M)\}.$ Hence, the result follows from \ref{corocotbf}.
\end{dem}
\vspace{.2cm}

We recall the following well known notions that will be useful in what follows.

\begin{defi} Let $\omega$ be a class of objects of the triangulated category $\T,$ and let ${\Omega}:=\cup_{i\in\Bbb{Z}}\,\omega[i].$  It is 
said that $\omega$ is a {\bf{cogenerator}} in $\T,$ if ${}^\perp{\Omega}=\{0\}.$ Dually, $\omega$ is a {\bf{generator}} in $\T,$ if 
${\Omega}^\perp=\{0\}.$
\end{defi}

\begin{rk}\label{gencogen}  Let $\omega$ be a class of objects of the triangulated category $\T.$ So, by induction and using the definition 
of $\overline{\Delta}_\T(\omega),$ it can be seen that $\omega$ is both a generator and a cogenerator in the triangulated category 
$\overline{\Delta}_\T(\omega).$
\end{rk}

\begin{pro}\label{weekcateg} Let $\X=\add\,(\X)$ be a cosuspended subcategory of $\T$ and let $\omega$ be an $\X$-injective weak-cogenerator in $\X.$  Then, $\cap_{i\in\Bbb{Z}}\,\X[i]=\{0\}$ if and only if  $\omega$ is  a cogenerator in
$\overline{\Delta}_\T(\X).$
\end{pro}
\begin{dem} First, by
\ref{M3(a)}, we have that $\overline{\Delta}_\T(\X)=\X^\wedge=\cup_{n\geq 0}\,\X[n]$ . We assert 
that $\cap_{i\in\Bbb{Z}}\,\X[i]\subseteq {}^\perp{\Omega}\cap \overline{\Delta}_\T(\X),$ where 
${\Omega}:=\cup_{i\in\Bbb{Z}}\,\omega[i].$ Indeed, let $M\in \cap_{i\in\Bbb{Z}}\,\X[i]$ and 
$j\in\Bbb{Z}.$ Hence $M=X[j-1]$ for some $X\in\X,$ and so $\Hom(M,W[j])\simeq\Hom(X,W[1])=0$ for 
any $W\in\omega,$ proving the assertion.
\

Assume that $\omega$ is  a cogenerator in $\overline{\Delta}_\T(\X).$ Hence ${}^\perp{\Omega}\cap \overline{\Delta}_\T(\X)=\{0\}$ and by the assertion above, it follows that $\cap_{i\in\Bbb{Z}}\,\X[i]=\{0\}.$
\

Suppose now  that  $\cap_{i\in\Bbb{Z}}\,\X[i]=\{0\}.$ Let $Y\in\overline{\Delta}_\T(\X)$ be non-zero. We prove the existence of an integer $\ell$ such that $\Hom(Y,\omega[\ell])\neq 0.$ Indeed, since $\overline{\Delta}_\T(\X)=\cup_{n\geq 0}\,\X[n],$ there is $n\in\Bbb{N}$ with $Y=X[n]$ for some $X\in \X.$ Furthermore, using that $X[n]=X[n-i][i]$ and the fact that $\X$ is cosuspended, it follows that $Y\in \X[i]$ for any $i\geq n.$ On the other hand, since $\cap_{j\in\Bbb{Z}}\,\X[j]=\{0\},$ we have that there is some $j_0<n$ such that $Y\notin \X[j_0].$ We assert that
$Y\notin \X[i]$ for any $i\leq j_0.$ It follows from $\X[i]=\X[i-j_0][j_0]\subseteq \X[j_0]$ and $Y\notin \X[j_0].$ Now, we set $\ell:=\mini\,\{s\;:\; j_0 < s\leq n\text{ and } Y\in\X[s]\}.$ So we have $Y[-\ell]\in \X$ and then, by using that  $\omega$ is a weak-cogenerator in $\X,$ there exists a distinguished triangle $X'[-1]\rightarrow Y[-\ell]\stackrel{f}{\rightarrow} W \rightarrow X'$ with $X'\in\X$ and $W\in\omega.$ Hence, the morphism
  $f:Y[-\ell]\rightarrow W$ is non-zero. In fact if $f=0,$ then  $Y[-\ell]$ would be a direct summand of $X'[-1]\in \X[-1],$ and so $Y[-\ell+1] \in  \X;$ giving a contradiction since $Y\notin \X[\ell-1].$ Thus $\Hom(Y,W[\ell])\neq0;$ proving the result.
\end{dem}

\begin{defi} Let $(\A, \B)$ be a co-$t$-structure on  $\T.$ It said that the pair $(\A, \B)$ is 
{\bf{non-degenerate below}} (respectively, {\bf{non-degenerate above}})  if \\ 
$\;\cap_{i\in\Bbb{Z}}\,\A[i]=\{0\}$ (respectively, $\;\cap_{i\in\Bbb{Z}}\,\B[i]=\{0\}$). So, 
it is said that $(\A, \B)$ is {\bf{non-degenerate}} if it is both non-degenerate below and above.
\end{defi}

\begin{pro}\label{critcog} Let $(\X,\Y)$ be a bounded below co-$t$-structure on a triangulated category $\T,$ and let $\omega:=\X\cap\Y.$ Then, the following conditions are equivalent.
\begin{itemize}
 \item[(a)] $(\X,\Y)$ is non-degenerate below.
 \item[(b)] $\omega$ is a cogenerator in $\T.$
\end{itemize}  
\end{pro}
\begin{dem} It follows from \ref{weekcateg}, \ref{cogenXiy} (a) and \ref{M3(a)} (c).
\end{dem}

\begin{cor}\label{corocritcog} Let $(\X,\Y)$ be a bounded co-$t$-structure on a triangulated category $\T,$ and let $\omega:=\X\cap\Y.$ Then, the following conditions are equivalent.
\begin{itemize}
 \item[(a)] $(\X,\Y)$ is non-degenerate.
 \item[(b)] $\omega$ is both a generator and a cogenerator in $\T.$
\end{itemize}
\end{cor}
\begin{dem} It follows from \ref{critcog} and its dual.
\end{dem}

\begin{cor}\label{ndcritcog} There is a bijective correspondence $\X\mapsto (\X,\overline{\Delta}_\T(\X)\cap\X^\perp[-1])$ between cosuspended subcategories $\X=\add\,(\X)$ of  $\T$ such that
$\X\cap\X^{\perp}[-1]$  is both a weak-cogenerator in $\X$ and a cogenerator in $\overline{\Delta}_\T(\X),$  and  non-degenerate below  co-$t$-structures $(\X,\Y)$  on   $\overline{\Delta}_\T(\X).$
\end{cor}
\begin{dem} From \ref{wcgext}, co-$t$-structures $(\X,\Y)$ on $\overline{\Delta}_\T(\X)$ correspond bijectively to cosuspended subcategories $\X$ of  $\T$ such that
$\X\cap\X^{\perp}[-1]$  is a weak-cogenerator in $\X.$ Therefore, the result follows from \ref{weekcateg} and \ref{critcog}.
\end{dem}
\vspace{.2cm}

The relationship between the different types of co-$t$-structures is as follows.

\begin{teo}\label{equivTipos}  Let $(\X,\Y)$ be a co-$t$-structure on a triangulated category $\T.$ Then, the following statements are equivalent.
\begin{itemize}
 \item[(a)] $(\X,\Y)$ is bounded.
 \item[(b)] $(\X,\Y)$ is faithful.
 \item[(c)] $(\X,\Y)$ is bounded and non-degenerate.
 \item[(d)] $\T=\overline{\Delta}_\T(\X\cap\Y).$
\end{itemize}
\end{teo}
\begin{dem} (a) $\Leftrightarrow$ (b) It is \ref{b=f}.
\

(a) $\Rightarrow$ (c) Assume that $(\X,\Y)$ is bounded. Thus, by \ref{bog} and 
\ref{criteriobf} (b), we get that $\T=\overline{\Delta}_\T(\omega)$ for $\omega:=\X\cap\Y;$ and so
by \ref{gencogen}, we have that $\omega$ is both a cogenerator and a generator in $\T.$ Then, (c) follows from \ref{corocritcog}.
\

(c) $\Rightarrow$ (d) It follows from \ref{corocotbf} (a).
\

(d) $\Rightarrow$ (a) Let $\T=\overline{\Delta}_\T(\X\cap\Y).$ Hence we get
$\T=\overline{\Delta}_\T(\X)=\overline{\Delta}_\T(\Y).$ Therefore, by \ref{M3(a)} and \ref{M3(a)dual} we get that $(\X,\Y)$ is bounded.
\end{dem}

\section{siltings and co-$t$-structures}

In this section, we show that in many cases a co-$t$-structure can be determined by a silting set. We also study the relationship between co-$t$-structures, silting and relative injective classes.
Following \cite{KV}, we recall the notion of a silting class in triangulated categories.

\begin{defi} Let $\omega$ be a class of objects in $\T.$ It is said that $\omega$  is {\bf{silting}} if $\id_{\omega}(\omega)=0.$
\end{defi}

We denote by ${}_\omega\sU$ (respectively, $\sU_\omega$) the smallest cosuspended (respectively, suspended) subcategory of $\T,$ closed under direct summands and containing $\omega.$

\begin{rk}\label{cosuspbar} Let $\omega$ be a class of objects in $\T.$ We define a sequence $\{\varepsilon^-_{i}(\omega)\}_{i\geq 0}$ of
classes of objects of  $\mathcal{T }$ as follows. Set $\varepsilon^-_{0}(\omega):=\add\,(\cup_{i\leq 0
 }\,\omega[i]).$ Assume that
$\varepsilon^-_{0}(\omega),\varepsilon^-_{1}(\omega),\cdots ,\varepsilon^-_{i-1}(\omega)$ are already defined. Then, we define $\varepsilon^-_{i}(\omega)$ as the class of objects in $\T,$ which are direct summands  of objects in  $\varepsilon^-_{i-1}(\omega)*\varepsilon^-_{0}(\omega).$ It is not hard to show that ${}_\omega\,\sU=\cup _{i\geq 0} \varepsilon^-_{i}(\omega).$
\end{rk}

\begin{lem}\label{Xomegads} Let $(\X,\omega)$ be a pair of classes of objects in $\T,$ such that $\omega\subseteq\X.$ Then, the following statements hold.
\begin{itemize}
 \item[(a)] If $\X$ is cosuspended and $\X=\add\,(\X),$  then $\X[-1]*\omega$ is closed under direct summands.
 \item[(b)] If $\omega$ is silting and closed under direct sums,  then $\omega$ is closed under extensions.
\end{itemize}
\end{lem}
\begin{dem} (a) Assume that $\X$ is cosuspended and closed under direct summands. Let $C\in \X[-1]*\omega.$ Then, there is a distinguished triangle $X[-1]\rightarrow
C\stackrel{f} \rightarrow W \rightarrow X$ where $X\in\X$ and $W\in\omega.$ Let $Z$ be a direct summand of $C,$ hence  there is distinguished triangle $Z\stackrel{u} \rightarrow C \rightarrow Z' \rightarrow Z[1],$ which splits.  Using the octahedral axiom, we get distinguished triangles $\Delta_1:\; Z\stackrel{fu} \rightarrow W \rightarrow V\rightarrow Z[1]$ and
 $\Delta_2:\;Z'\rightarrow  V \rightarrow X\rightarrow Z'[1].$ By the hypothesis, we have that $\X[-1]*\omega\subseteq\X*\X\subseteq\X;$ and so $C\in\X,$ giving us that $Z$ and $Z'$ belong to $\X.$ Thus $V\in\X$ (see $\Delta_2$), and hence from $\Delta_1,$ we get that $Z\in\X[-1]*\omega.$
\

(b) Assume that $\omega$ is silting and closed under direct sums. Let $\Delta:\; W\to X\to W'\to W[1]$ be a distinguished triangle with $W,W'\in\omega.$ Using that $\id_\omega(\omega)=0,$ we obtain that the triangle $\Delta$ splits; and hence $X\in\omega$ since $\omega$ is closed under direct sums.
\end{dem}

\begin{pro}\label{siltingcog} Let $\omega$ be a silting class in $\T$ such that $\add\,(\omega)=\omega.$ Then $\omega$ is an ${}_\omega\,\sU$-injective weak-cogenerator in  ${}_\omega\,\sU.$
\end{pro}
\begin{dem} From \ref{cosuspbar}, we know that  ${}_\omega\,\sU=\cup _{n\geq 0} \varepsilon^-_{n}(\omega).$  Hence, it is enough to prove, by induction on $n,$ that $\varepsilon^-_{n}(\omega) \subseteq {}_\omega\, \sU[-1]*\omega$ for any $n\in\Bbb{N}.$ Assume that $\add\,(\omega)=\omega.$ In particular, we have that $\varepsilon^-_{0}(\omega)=\oplus_{i\leq 0}\,\omega [i],$ where direct sums means here finite direct sums.

 If $X\in \varepsilon^-_{0}(\omega),$ then there is an split distinguished triangle $W'\rightarrow X \rightarrow W\rightarrow W'[1],$     where $W'\in \oplus _{i< 0}\,\omega[i]$ and $W\in\omega.$ Hence $X\in {}_\omega\,\sU[-1]*\omega.$

Let $n>1,$ and take $X\in \varepsilon^-_{n}(\omega).$ Then, there is a distinguished triangle $ X_{n-1}\rightarrow
X^{\prime }\rightarrow X_{0} \rightarrow X_{n-1}[1]$ with $X_{0}\in\varepsilon^-_{0}(\omega),$ $X_{n-1}\in\varepsilon^-_{n-1}(\omega)$ and $X$ is a direct summand of $X'.$  For $X_0$ we have an  split distinguished triangle $W'\rightarrow X_0 \stackrel{f}{\rightarrow} W\rightarrow W'[1],$     where $W'\in \oplus _{i< 0}\,\omega[i]$ and $W\in\omega.$ Therefore, by the base change argument (using the octahedral axiom), we get the following commutative and exact diagram in $\T$

$$\begin{CD}
@.   W[-1] @= W[-1] @.\\
@. @VVV @VVV @.\\
X_{n-1}@>>> Y @>>>  W'  @>>> X_{n-1} [1]\\
@| @VVV @VVV @|\\
X_{n-1} @>>> X' @>>>  X_0 @>>> X_{n-1} [1]\\
@. g@VVV @VVfV @.\\
@.  W @=  W  @.\,  \\
{}
  \end{CD}
  $$
By induction there exist a distinguished triangle  $U[-1] \rightarrow X_{n-1} \stackrel{h}\rightarrow W''\rightarrow U$ where $U\in{}_\omega\,\sU$ and $W''\in\omega.$ Since $\Hom(\oplus_{i< 0}\,\omega [i],\omega[1])=0$ because $\omega$ is silting, we have a morphism $\alpha:W'\to U$ that can be completed to a distinguished triangle $W'\stackrel{\alpha}{\rightarrow}U\stackrel{\beta}{\rightarrow}V\rightarrow W'[1].$ By using the octahedral axiom, we get the following exact and commutative diagram in $\T$

$$\begin{CD}
U[-1] @>\beta[-1]>>  V[-1] @>>>  W'  @>\alpha>>  U\\
@VVV @VVV @| @VVV  \\
X_{n-1}@>>> Y @>>>  W'  @>>> X_{n-1} [1]\\
@VhVV @VVV @VVV @VVh[1]V\\
 W'' @=  W''  @>>> 0 @>>> W''[1]\\
@VVV @VVV @.\\
U  @>>> V @.\\
{}
  \end{CD}$$

From the triangle $U[-1] \rightarrow V[-1]\rightarrow  W'\stackrel{\alpha}{\rightarrow} U,$ it follows that $V[-1]\in{}_\omega\,\sU[-1]$ since $ {}_\omega\,\sU[-1]$ is closed under extensions. Now,  the triangle $V[-1]\rightarrow Y \stackrel{} \rightarrow W''\rightarrow V$ implies that $Y \in {}_\omega\,\sU[-1]*\omega .$ Then $X' \in  {}_w \U[-1] *\omega * \omega $ since we have the triangle $W[-1] \rightarrow Y\rightarrow  X'\stackrel{g}{\rightarrow} W.$ But  $ \omega*\omega \subseteq\omega$ (see \ref{Xomegads} (b)), and so  $X'\in{}_\omega\,\sU[-1] * \omega * \omega \subseteq{}_\omega\,\sU[-1]*\omega .$ Therefore, from \ref{Xomegads} (a), we conclude that $X\in {}_\omega\,\sU[-1]*\omega;$ proving that $\omega$ is a weak-cogenerator in  ${}_\omega\,\sU.$ Finally, we prove that $\omega$ is also ${}_\omega\,\sU$-injective. Indeed, since $\id_\omega(\omega)=0$ it follows from \cite[Lemma 4.2 (a2)]{MSSS1} that $\omega\subseteq {}_\omega\,\U^\perp[-1];$ and using that 
${}_\omega\,\U^\perp[-1]={}_{{}_\omega\,\sU}\,\U^\perp[-1],$ we get by \cite[Lemma 4.2 (a2)]{MSSS1} that $\id_{{}_\omega\,\sU}\,(\omega)=0.$
\end{dem}

The following result is very similar to \cite[Theorem 4.3.2(II)]{Bo}, which was
proved with different techniques.

\begin{teo}\label{siltingteo} Let $\omega$ be a silting class in $\T$ such that $\omega=\add\,(\omega).$ Then, $\omega={}_\omega\,\U\cap \U_{\,\omega}$ and the pair 
$({}_\omega\,\U, \U_{\,\omega})$ is a bounded co-t-structure on $\thick(w).$  
\end{teo}
\begin{dem} Since $\thick(\omega)=\thick({}_\omega\,\sU),$ it follows from \ref{M3(a)} that  $\thick(\omega)={}_{\omega}\,\sU^\wedge.$ On the other 
hand, by \ref{siltingcog} and \ref{correspond} (a), we get that the pair $({}_{\omega}\,\sU,\omega^\wedge)$ is a co-t-structure on 
$\thick({}_\omega\,\sU)=\thick(\omega)$ and $\omega={}_{\omega}\,\sU\cap\omega^\wedge.$ In particular, from \ref{criterio} (a), it follows that 
$\U_{\,\omega}=\omega^\wedge$ and hence $\U_{\,\omega}=\sU_{\,\omega}.$ Therefore, the pair $({}_{\omega}\,\sU, \sU_{\,\omega})$ is a bounded 
below and faithful below co-t-structure on $\thick(w),$ and $\omega={}_{\omega}\,\sU\cap \sU_{\,\omega}.$ So, from \ref{criteriobf}, we get 
that $({}_{\omega}\,\sU, \sU_{\,\omega})$ is bounded on $\thick(w).$ Furthermore, by \ref{corocotbf} (a), we obtain that 
${}_\omega\,\U={}_\omega\,\sU$ and $\U_{\,\omega}=\sU_{\,\omega}.$  
\end{dem}

\begin{rk}\label{rksiltingteo} Let $\omega$ be a silting class in $\T$ such that $\omega=\add\,(\omega).$ Then, by \ref{siltingteo}, it follows that ${}_\omega\,\U={}_\omega\,\sU$ and $\U_{\,\omega}=\sU_{\,\omega}.$
\end{rk}

\begin{defi} For a given triangulated category $\T,$ we introduce the following classes:
\begin{itemize}
 \item[(a)] $\bf{S}$ consists of all silting classes $\omega$ of $\T$ such that $\add\,(\omega)=\omega.$
 \item[(b)] $\C_b$ consists of all bounded co-$t$-structures $(\X,\Y)$ on $\thick(\X\cap\Y).$
 \end{itemize}
\end{defi}

\begin{cor}\label{siltingcorresp} Let $\T$ be a triangulated category. Then, the correspondence $\varphi:\bf{S}\to\C_b,$ given by $\varphi(\omega):=({}_\omega\,\U, \U_{\,\omega}),$ is bijective.
\end{cor}
\begin{dem} From \ref{siltingteo}, it follows that $\varphi:\bf{S}\to\C_b$ is well defined and injective. Let $(\X,\Y)$ in $\C_b,$ and consider $\omega:=\X\cap\Y.$ Since $(\X,\Y)$ is a bounded  co-$t$-structure on $\thick(\omega),$ we conclude by \ref{corocotbf} (a) that $\varphi(\omega)=(\X,\Y);$ proving that $\varphi$ is also surjective.
\end{dem}

\begin{cor}\label{correspK} Let $\T$ be a triangulated category. Then, there is a bijective
correspondence $(\X,\Y)\mapsto\omega:=\X\cap\Y,$ with inverse 
$\omega\mapsto ({}_\omega\,\U,\U_{\,\omega}),$
between bounded co-$t$-structures $(\X,\Y)$ on $\T$ and silting classes $\omega=\add\,(\omega)$  such that $\T=\thick(\omega).$
\end{cor}
\begin{dem} It follows from \ref{siltingcorresp} and \ref{equivTipos}.
\end{dem}
\vspace{.2cm}

The next result characterizes when a cosuspended subcategory of $\T$ determines a bounded co-$t$-structure on $\T.$

\begin{teo}\label{Msc} Let $\T$ be a triangulated category, and $\X$ be a cosuspended subcategory of $\T$ such that $\X=\add\,(\X).$ Then, the following statements are equivalent.
\begin{itemize}
 \item[(a)] There is a bounded co-$t$-structure $(\X,\Y)$ on $\T.$
 \item[(b)]  $\thick(\X\cap\X^\perp[-1])=\T.$
 \item[(c)] There is an $\X$-injective $\omega=\add\,(\omega)$ such that $\thick(\omega)=\T$  and $\omega\subseteq\X.$
 \item[(d)] There is a silting $\omega=\add\,(\omega)$ such that $\thick(\omega)=\T\;$ and 
 $\;\omega\subseteq\X\subseteq\omega^\vee.$
\end{itemize}
Moreover, if one of the above conditions hold, we have that $\X={}_\omega\,\U=\omega^\vee,$ $\Y=\U_{\,\omega}\;$ and $\;\omega=\X\cap\Y=\X\cap\X^\perp[-1].$
\end{teo}
\begin{dem} (a) $\Rightarrow$ (d) Assume that $(\X,\Y)$ is a bounded co-$t$-structure on $\T,$ 
and let $\omega=\X\cap\Y.$ Then, by \ref{correspK}, we get that $\omega=\add\,(\omega)$ is a silting  class such that $\thick(\omega)=\T.$ On the other hand, since $(\X,\Y)$ is bounded, it follows from \ref{corocotbf} (a) that $\X={}_\omega\,\U=\omega^\vee$ and $\Y=\U_{\,\omega}.$ 
\

(d) $\Rightarrow$ (a) Suppose there is a silting class $\omega$ such that $\omega\subseteq\X\subseteq\omega^\vee$ and
$\thick(\omega)=\T.$
Hence, by \ref{siltingteo}, it follows that $({}_\omega\,\U,\U_{\,\omega})$ is a bounded co-$t$-structure on $\T,$ and $\omega={}_\omega\,\U\cap\U_{\,\omega}.$ In particular, from \ref{corocotbf} (a), we know that ${}_\omega\,\U=\omega^\vee.$ Furthermore, since $\omega\subseteq\X,$ it follows that ${}_\omega\,\U\subseteq\X$ and hence $\X={}_\omega\,\U.$ 
\

(a) $\Rightarrow$ (b) Let $(\X,\Y)$ be bounded. Then, by \ref{corocotbf} (b), we get that 
 $\thick(\omega)=\T,$ where $\omega:=\X\cap\Y=\X\cap \X^\perp[-1]$ (see \ref{proP} (c)). 
\

(b) $\Rightarrow$ (c) Let $\thick(\X\cap\X^\perp[-1])=\T.$ Since $\X$ is a cosuspended subcategory of $\T,$ it follows from \cite[Lemma 4.2 (a2)]{MSSS1}, that $\omega:=\X\cap\X^\perp[-1]$ is $\X$-injective.
\

(c) $\Rightarrow$ (a) Assume the hypothesis in (c). In particular, $\omega$ is silting since 
$\id_\X(\omega)=0.$ Thus, 
from \ref{siltingteo}, it follows that 
$({}_\omega\,\U,\U_{\,\omega})$ is a bounded co-$t$-structure on $\T$ and also 
that $\omega={}_\omega\,\U\cap\U_{\,\omega}.$ Furthermore ${}_\omega\,\U\subseteq\X$ since 
$\omega\subseteq\X.$ On the other hand, since $\pd_\omega(\X)=\id_\X(\omega)=0,$  it follows 
from \cite[Lemma 4.2 (a1)]{MSSS1}, that $\X\subseteq {}^\perp\U_{\,\omega}[1]={}_\omega\,\U$ (see
\ref{proP} (c)); and 
hence $\X={}_\omega\,\U.$ 
\end{dem}

\section{co-$t$-structures on $\D(\HH)$}

Throughout this section, $k$ denotes an algebraically closed field and $\HH$ an {\emph{abelian 
hereditary $k$-category which is Hom-finite, Ext-finite and has a tilting object}}. We will 
consider the bounded derived category $\D(\HH)$ which is triangulated and has been intensively 
studied (see, for example, \cite{H2} and \cite{HR}). 
\

In this section, we give a description of the bounded co-$t$-structures on 
$\T:=\D(\HH).$ In this case, the obtained results take a more complete form that in the preceding section.
\

In what follows, we need the following useful lemma. For details, we refer the reader to \cite{AST}.

\begin{lem}\cite{AST} \label{heredgen} Let $\omega$ be a set in the triangulated category $\D(\HH).$ Then, the following statements hold.
\begin{itemize}
 \item[(a)] $\omega$ is a generator in $\D(\HH)$ if and only if it is a cogenerator in $\D(\HH).$
 \item[(b)] Let $\omega$ be a silting  class in $\D(\HH).$ Then, $\omega$ is a generator in $\D(\HH)$ if and only if $\thickH(\omega)=\D(\HH).$
 \end{itemize}
\end{lem}
\begin{dem} (a) It follows from \cite[Lemma 2.1]{AST} since $\D(\HH)$ has a Serre duality.
\

(b) ($\Leftarrow$) It follows from \ref{gencogen}.
\

($\Rightarrow$) By \cite[Corollary 3.2 (b)]{AST}, we know that, for every complex $X\in \T,$ there is a distinguished triangle $W\rightarrow X \rightarrow L \rightarrow W[1]$ such that $W\in \thickH(\omega)$ and $L\in \thickH(\omega)^{\perp}.$
If $\omega$ is a generator then $L=0$ and so $X\simeq W\in \thick(\omega),$ proving that $\thickH(\omega)=\D(\HH).$
\end{dem}

\begin{pro}\label{equivBnondegH} Let $(\X,\Y)$ be a co-$t$-structure on $\D(\HH).$ Then, the 
following statements are equivalent.
\begin{itemize}
 \item[(a)] $(\X,\Y)$ is bounded.
 \item[(b)] $(\X,\Y)$ is non-degenerate below and bounded below.
 \item[(c)] $(\X,\Y)$ is non-degenerate above and bounded above.
 \end{itemize}
\end{pro}
\begin{dem} (a) $\Rightarrow$ (b) It follows from \ref{equivTipos}.
\

(b) $\Rightarrow$ (a) Assume the hypothesis in (b). Then, by \ref{critcog}, we get that $\omega$ 
is a cogenerator in $\D(\HH).$ Hence, from \ref{heredgen} and \ref{equivTipos}, we conclude 
that $(\X,\Y)$ is bounded.
\

Finally, the equivalence between (c) and (a), can be proven in a similar way we did for (a) and (b).
\end{dem}
\vspace{.2cm}

The following result gives a characterization, in terms of generators and cogenerators, of the 
bounded co-$t$-structures on $\D(\HH).$

\begin{teo}\label{MscH} Let $\X$ be a cosuspended subcategory of $\D(\HH)$ such that $\X=\add\,(\X).$ Then, the following statements are equivalent.
\begin{itemize}
 \item[(a)] There is a bounded co-$t$-structure $(\X,\Y)$ on $\D(\HH).$
 \item[(b)]  $\X\cap\X^\perp[-1]$ is a generator set in $\D(\HH).$
 \item[(c)] There is an $\X$-injective set $\omega=\add\,(\omega),$ which is a cogenerator in $\D(\HH)$ and $\omega\subseteq\X.$
\end{itemize}
\end{teo}
\begin{dem} Since $\X$ is cosuspended, it follows by \cite[Lemma 4.2 (a2)]{MSSS1} that $\X\cap\X^\perp[-1]$ is $\X$-injective; and hence, it is silting. Therefore, the result follows from \ref{Msc} and \ref{heredgen}.
\end{dem}

\begin{cor}\label{precovMsc} Let $\X$ be a cosuspended subcategory of $\D(\HH)$ such that $\X=\add\,(\X).$ If $\X\cap\X^\perp[-1]$ is 
a generator set in $\T,$ then $\X^\wedge=\T$ and $\X$ is a precovering class in $\T.$
\end{cor}
\begin{dem} It follows from \ref{MscH} (a) and \ref{correspond} (b).
\end{dem}

\begin{cor}\label{coroMsc} Let $(\X,\omega)$ be a pair of classes of objects of $\D(\HH),$ which are closed under direct summands, and let
$\X$ be cosuspended. Then, the following conditions are equivalent.
\begin{itemize}
 \item[(a)] $\omega$ is an $\X$-injective weak-cogenerator in $\X,$ $\bigcap_{i\in\Enteros}\,\X[i]=\{0\}$ and $\X^\wedge=\D(\HH).$
 \item[(b)] $\omega\subseteq\X\subseteq\omega^\vee\,$ and $\,\omega=\add\,(\omega)$  is a silting cogenerator set in $\D(\HH).$
 \item[(c)] $\omega=\add\,(\omega)\subseteq\X$ and $\omega$ is an $\X$-injective cogenerator set 
 in $\D(\HH).$
 \item[(d)] $\omega=\X\cap\X^\perp[-1]$ and $\omega$ is a generator set in $\D(\HH).$
\end{itemize}
Moreover, if one of the above conditions hold, we have that $\X={}_\omega\,\U=\omega^\vee.$
\end{cor}
\begin{dem} (a) $\Rightarrow$ (b) By \ref{correspond} (a) and \ref{M3(a)}, there is a bounded 
below co-$t$-structure $(\X,\Y)$ on $\D(\HH).$ Furthermore, by \ref{weekcateg}, we get that 
$\omega$ is a cogenerator set in $\D(\HH);$ and so, by \ref{heredgen} and \ref{equivTipos}, it 
follows that $(\X,\Y)$ is bounded. Hence (b) follows from \ref{Msc} (d).
\

(b) $\Rightarrow$ (a) By \ref{heredgen}, it follows that $\thickH(\omega)=\D(\HH).$ Therefore, the
condition (d) in \ref{Msc} holds. Hence (a) follows by \ref{correspond} (a) and \ref{weekcateg}.
\

(b) $\Leftrightarrow$ (c) It follows from \ref{heredgen}, \ref{Msc} and \ref{MscH}.
\

(a) $\Leftrightarrow$ (d) It follows from \ref{MscH}, \ref{cogenXiy} (a), \ref{heredgen} and \ref{weekcateg}.
\end{dem}
\vspace{.2cm}

Let $\omega$ be a class of objects of $\D(\HH).$ We say that $\omega$ is of {\bf{finite type}} if there exist a finite number of pairwise non isomorphic indecomposable objects $W_1,W_2,\cdots,W_n$ in $\D(\HH)$ satisfying that $\add\,(\omega)=\add\,(\{W_1,W_2,\cdots,W_n\}).$ In such a case, we set $\mathrm{ind}\,(\omega):=\{W_1,W_2,\cdots,W_n\}$ and $\mathrm{rk}\,(\omega):=n.$ We also denote by $\mathrm{rk}\,K_0\,(\HH)$ 
the rank of the Grothendieck group associated with $\HH.$

\begin{lem}\label{auxteoH}\cite{AST} The following statements holds.
 \begin{itemize}
  \item[(a)] If $\omega$ is a silting set in $\D(\HH),$ then $\mathrm{rk}\,(\omega)\leq \mathrm{rk}\,K_0\,(\HH).$
  \item[(b)] Let $\Y=\add\,(\Y)$ be a suspended and precovering subcategory of $\D(\HH),$ and let $\omega:=\Y\cap{}^\perp\Y[1].$ Then, $\mathrm{rk}\,(\omega)=\mathrm{rk}\,K_0\,(\HH)$ if and only if $\omega$ is a generator in $\D(\HH).$ Furthermore, if this is the case, then $\Y=\sU_{\,\omega}=\U_{\,\omega}.$
 \end{itemize}
\end{lem}
\begin{dem} (a) By \ref{siltingcog}, we know that $\omega$ is ${}_\omega\,\sU$-injective. 
Therefore, the item (a) is just the dual of \cite[Theorem 2.3 (b)]{AST}.
\

(b) This is \cite[Corollary 4.4]{AST}. Observe that the equality $\sU_{\,\omega}=\U_{\,\omega}$ follows from \ref{rksiltingteo}
\end{dem}

\begin{teo}\label{teoH} There are bijective correspondences $$(\X,\Y)\mapsto \Y,\quad \Y\mapsto\omega:=\Y\cap{}^\perp\Y[1]\quad\text{ and }\quad\omega\mapsto({}_\omega\,\U,\U_{\,\omega})$$  between the following classes:
 \begin{itemize}
  \item[(a)] Bounded co-$t$-structures $(\X,\Y)$ on $\D(\HH).$
  \item[(b)] Suspended and precovering subcategories $\Y=\add\,(\Y)$ of $\D(\HH)$ such that  $\mathrm{rk}\,(\Y\cap{}^\perp\Y[1])=\mathrm{rk}\,K_0\,(\HH).$
  \item[(c)] Silting sets $\omega=\add\,(\omega)$ in $\D(\HH)$ such that $\mathrm{rk}\,(\omega)=\mathrm{rk}\,K_0\,(\HH).$
 \end{itemize}
\end{teo}
\begin{dem} By \cite[Corollary 4.5]{AST} and \ref{auxteoH} (b), we have that the correspondence $\Y\mapsto\Y\cap{}^\perp\Y[1]$ between the classes of items (b) and (c) is bijective with inverse $\omega\mapsto\U_{\,\omega}.$
\

We prove now that the correspondence $(\X,\Y)\stackrel{\alpha}{\mapsto}\X\cap\Y$ between the classes of items (a) and (c) is bijective with inverse $\omega\stackrel{\beta}{\mapsto}({}_\omega\,\U,\U_{\,\omega}).$ Indeed, let $(\X,\Y)$ be a pair belonging to item (a). By \ref{Msc} and \ref{MscH}, it follows that $\X\cap\Y$ is a silting generator in $\D(\HH)$ and $\Y=\U_{\X\cap\Y}.$ Hence, by applying \cite[Corollary 3.2 (b)]{AST}, we get that $\Y$ is a suspended and precovering subcategory of $\D(\HH).$ Therefore, from
\ref{auxteoH}, we get that $\X\cap\Y$ belongs to the item (c). Furthermore, from \ref{Msc}, we conclude that $\beta\,\alpha(\X,\Y)=(\X,\Y).$ Let $\omega$ be a class belonging to the item (c). In particular, by \ref{siltingteo}, we have that $\beta(\omega)=({}_\omega\,\U,\U_{\,\omega})$ is a bounded non-degenerate co-$t$-structure on $\thick(\omega)$ and $\omega={}_\omega\,\U\cap\U_{\,\omega}=\alpha\,\beta(\omega).$ But, using the bijective correspondence between the classes of items (b) and (c), we get that $\U_{\,\omega}$ is a suspended and precovering subcategory of $\D(\HH).$ Therefore, from
\ref{auxteoH}, we obtain that $\omega$ is a generator in $\D(\HH);$ and so $\thick(\omega)=\D(\HH)$ (see \ref{heredgen}), proving that $({}_\omega\,\U,\U_{\,\omega})$ is a bounded co-$t$-structure on $\D(\HH).$ That is, $\beta(\omega)$ belongs to the item (a).
\end{dem}

\begin{rk}\label{rkteoH} The item $(\mathrm{b})$ in \ref{teoH} is equivalent to the following one:
\

$(\mathrm{b'})$ Suspended subcategories $\Y=\add\,(\Y)$ of $\D(\HH)$ such that  $\mathrm{rk}\,(\Y\cap{}^\perp\Y[1])=\mathrm{rk}\,K_0\,(\HH).$
\

\noindent Moreover, if $(\mathrm{b'})$ holds, then we have that $\omega:=\Y\cap{}^\perp\Y[1]$ is a generator set in $\D(\HH)$ and $\Y=\U_\omega.$
\end{rk}
\begin{dem} Let $\Y=\add\,(\Y)$ be a suspended subcategory of $\D(\HH),$ and let $\omega:=\Y\cap{}^\perp\Y[1]$ be such that  $\mathrm{rk}\,(\omega)=\mathrm{rk}\,K_0\,(\HH).$ Then, from \ref{teoH} (a), we have that $({}_\omega\,\U,\U_{\,\omega})$ is a bounded co-$t$-structure on $\D(\HH).$ Thus $\overline{\Delta}_{\D(\HH)}(\omega)=\D(\HH)$ and $\omega$ is a generator set in $\D(\HH)$ (see \ref{Msc}). In particular $(\overline{\Delta}_{\D(\HH)}(\omega))^\perp=\{0\};$ and therefore, from \cite[Theorem 4.2 (b)]{AST}, we conclude that $\Y=\U_\omega.$ Finally, using \cite[Corollary 3.2]{AST}, we get that $\Y$ is precovering in $\D(\HH).$
\end{dem}

\begin{cor}\label{c1teoH} There are bijective correspondences $$(\X,\Y)\mapsto \X,\quad \X\mapsto\omega:=\X\cap\X^\perp[-1]\quad\text{ and }\quad\omega\mapsto({}_\omega\,\U,\U_{\,\omega})$$ between the following classes:
 \begin{itemize}
  \item[(a)] Bounded co-$t$-structures $(\X,\Y)$ on $\D(\HH).$
  \item[(b)] Cosuspended and preenveloping subcategories $\X=\add\,(\X)$ of $\D(\HH)$ such that  $\mathrm{rk}\,(\X\cap\X^\perp[-1])=\mathrm{rk}\,K_0\,(\HH).$
  \item[(c)] Silting sets $\omega=\add\,(\omega)$ in $\D(\HH)$ such that $\mathrm{rk}\,(\omega)=\mathrm{rk}\,K_0\,(\HH).$
  \item[(d)] Cosuspended subcategories $\X=\add\,(\X)$ of $\D(\HH)$ such that  $\mathrm{rk}\,(\X\cap\X^\perp[-1])=\mathrm{rk}\,K_0\,(\HH).$
 \end{itemize}
\end{cor}
\begin{dem} Let $\T:=\D(\HH).$ In order to prove the result, using \ref{teoH}, \ref{rkteoH} and the duality principle for triangulated categories, it is enough to prove the following statement: if $(\X,\Y)$ is a bounded co-$t$-structure on $\T,$ then $(\Y^{op},\X^{op})$ is so on the opposite triangulated category $\T^{op}.$ Observe, firstly, that this statement is true since the boundedness property is a self-dual notion; and secondly,  
$\T^{op}\simeq\D(\HH^{op})$ where $\HH^{op}$ is also  an abelian 
hereditary $k$-category which is Hom-finite, Ext-finite and has a tilting object (see \cite[Proposition 1.9]{HR}).
\end{dem}
\vspace{.2cm}

As a nice consequence, from \ref{teoH} (b) and \ref{c1teoH} (b), is the following corollary, 
saying that any bounded co-$t$-structure on $\D(\HH)$ has two companions as 
$t$-structures: one on the left and the other on the right. For the convenience of the reader, we recall the definition of $t$-structure.

\begin{defi}\cite{BBD} A pair $(\A,\B)$ of subcategories in $\T$ is said to be a {\bf{$t$-structure}} on $\T$ if the following conditions hold.
 \begin{itemize}
  \item[(a)] $\A[1]\subseteq\A$ and $\B[-1]\subseteq \B.$
  \item[(c)] $\Hom_\T(\A,\B[-1])=0.$
  \item[(d)] $\T=\A*\B[-1].$
 \end{itemize}
\end{defi}

\begin{cor}\label{adjacent} Let $(\X,\Y)$ be a bounded co-$t$-structure on $\D(\HH).$ Then, the pairs $({}^\perp\X[-1],\X)$ and $(\Y,\Y^\perp[1])$ are both $t$-structures on  $\D(\HH).$
\end{cor}
\begin{dem} From \ref{teoH} (b) and \cite[Proposition 1.3]{KV}, it follows that $\Y$ is an aisle 
in $\D(\HH).$ Thus $(\Y,\Y^\perp[1])$ is a $t$-structure on  $\D(\HH).$ Furthermore, by \ref{c1teoH} (b) and the dual of \cite[Proposition 1.3]{KV},  we get that $\X$ is a co-aisle 
in $\D(\HH),$ and so, $({}^\perp\X[-1],\X)$ is a $t$-structure on  $\D(\HH).$
\end{dem}

\begin{rk} The previous result, says that a bounded co-$t$-structure on $\D(\HH)$ is always left (respectively, right) adjacent to a $t$-structure on $\D(\HH)$ in the sense of \cite{Bo} 
\end{rk}

\begin{cor}\label{c2teoH} Let $(\X,\Y)$ be a bounded co-$t$-structure on $\D(\HH).$ Then $\X$ and $\Y$ are functorially finite in $\D(\HH).$ 
\end{cor}
\begin{dem} It follows from \ref{teoH}, \ref{c1teoH} and \ref{proP}.
\end{dem}

\begin{cor}\label{c3teoH} Let $\omega=\add\,(\omega)$ be a silting generator set in $\D(\HH).$ Then ${}_\omega\,\U$ and $\U_{\,\omega}$ are functorially finite in $\D(\HH),$  ${}_\omega\,\U^\wedge=\D(\HH)=\U_{\,\omega}^\vee$ and $\mathrm{rk}\,(\omega)=\mathrm{rk}\,K_0\,(\HH).$
\end{cor}
\begin{dem} From \ref{heredgen} (b) and  \ref{siltingteo}, we know that $({}_\omega\,\U,\U_{\,\omega})$ is a bounded co-$t$-structure on $\D(\HH).$  Hence the result follows from \ref{c2teoH} and \ref{teoH}
(c).
\end{dem}

\begin{cor}\label{c4teoH} Let $\omega$ be a silting set in $\D(\HH).$ Then,  $\omega$ is a generator in $\D(\HH)$ if and only if  $\mathrm{rk}\,(\omega)=\mathrm{rk}\,K_0\,(\HH).$
\end{cor}
\begin{dem} Consider $\omega':=\add\,(\omega).$ Observe that $\omega':=\add\,(\omega')$ and $\omega'$ is also a silting set in $\D(\HH).$
\

($\Rightarrow$) Suppose that $\omega$ is a generator in $\D(\HH);$ and hence $\omega'$ is so. Then by \ref{c3teoH}, we get $\mathrm{rk}\,(\omega)=\mathrm{rk}\,K_0\,(\HH)$ since $\mathrm{rk}\,(\omega)=\mathrm{rk}\,(\omega').$
\

($\Leftarrow$) Assume now that $\mathrm{rk}\,(\omega)=\mathrm{rk}\,K_0\,(\HH).$ Thus $\mathrm{rk}\,(\omega')=\mathrm{rk}\,K_0\,(\HH)$ and so from \ref{teoH}, it follows that $({}_{\omega'}\,\U,\U_{\,\omega'})$ is a bounded co-$t$-structure on $\D(\HH).$ Therefore by 
\ref{criteriobf} and \ref{b=f}, we get that $\overline{\Delta}_{\D(\HH)}(\omega)=\overline{\Delta}_{\D(\HH)}(\omega')=\D(\HH).$ Hence $\omega$ is a generator in $\D(\HH)$ (see \ref{heredgen}).
\end{dem}
\begin{cor}\label{c5teoH} Let $\Y=\add\,(\Y)$ be a suspended subcategory of $\D(\HH)$ and let $\omega:=\Y\cap{}^\perp\Y[1].$ If   $\mathrm{rk}\,(\omega)=\mathrm{rk}\,K_0\,(\HH),$ then $\Y$ is functorially finite in $\D(\HH),$ $\Y=\U_\omega$ and $\Y^\vee=\D(\HH).$
\end{cor}
\begin{dem} Let $\mathrm{rk}\,(\omega)=\mathrm{rk}\,K_0\,(\HH).$ Then, by \ref{rkteoH}, it follows that  $\omega$ is a generator set in $\D(\HH)$ and $\Y=\U_\omega.$ So the result now follows from \ref{c3teoH}.
\end{dem}

\begin{cor}\label{c65teoH} Let $\X=\add\,(\X)$ be a cosuspended subcategory of $\D(\HH)$ and let $\omega:=\X\cap\X^\perp[-1].$ If   
$\mathrm{rk}\,(\omega)=\mathrm{rk}\,K_0\,(\HH),$ then $\X$ is functorially finite in $\D(\HH),$ $\X={}_\omega\U$ and $\X^\wedge=\D(\HH).$
\end{cor}
\begin{dem} It follows from \ref{c5teoH} and the discussion given in the proof of \ref{c1teoH}. 
\end{dem}
\vspace{.4cm}

\textbf{Acknowledgement} The authors are very grateful to the referee for the comments, corrections and suggestions.

\footnotesize

\vskip3mm \noindent Octavio Mendoza Hern\'andez:\\ Instituto de Matem\'aticas, Universidad Nacional Aut\'onoma de M\'exico\\
Circuito Exterior, Ciudad Universitaria,
C.P. 04510, M\'exico, D.F. MEXICO.\\ {\tt omendoza@matem.unam.mx}

\vskip3mm \noindent Edith Corina S\'aenz Valadez:\\ Departamento de Matem\'aticas,
Facultad de Ciencias, Universidad Nacional Aut\'onoma de M\'exico\\ Circuito Exterior, Ciudad
Universitaria,
C.P. 04510, M\'exico, D.F. MEXICO.\\ {\tt ecsv@lya.fciencias.unam.mx}

\vskip3mm \noindent Valente Santiago Vargas:\\ Instituto de Matem\'aticas, Universidad Nacional Aut\'onoma de M\'exico\\
Circuito Exterior, Ciudad Universitaria,
C.P. 04510, M\'exico, D.F. MEXICO.\\ {\tt valente@matem.unam.mx}

\vskip3mm \noindent Mar\'{i}a Jos\'e Souto Salorio:\\ Facultade de Inform\'atica,
Universidade da Coru\~na\\ 15071 A Coru\~na, ESPA\~NA.\\ {\tt mariaj@udc.es}

\end{document}